\documentclass[12pt]{elsarticle}
\usepackage{amsfonts,amsmath,amssymb,amsthm,amscd,bbm}
\usepackage{listings}
\usepackage{graphicx}

\usepackage{natbib}{\tiny }
\usepackage[usenames,dvipsnames]{color}

\def\ignore#1{{}}
\def\st#1{\textcolor{red}{#1}}

\def\bbbe{\mathbb{E}}
\def\bbbp{\mathbb{P}}

\usepackage{bm}
\begin{document}

\begin{frontmatter}

\title{Ancestral inference from haplotypes and mutations} 



\author{Robert C. Griffiths}
\address{Department of Statistics, University of Oxford, 24--29 St Giles, Oxford OX1 3LB, UK}

\author{Simon Tavar\'e}
\address{DAMTP, University of Cambridge, Centre for Mathematical Sciences, Wilberforce Road, Cambridge CB3 0WA, UK}





\begin{abstract}
We consider inference about the  history of a sample of DNA sequences, conditional upon the haplotype counts and the number of segregating sites observed at the present time. After deriving some theoretical results in the coalescent setting, we implement rejection sampling and importance sampling schemes to perform the inference. The importance sampling scheme addresses an extension of the Ewens Sampling Formula for a configuration of haplotypes and the number of segregating sites in the sample. The implementations include both constant and variable population size models. The methods are illustrated by two human Y chromosome data sets.
\end{abstract}


\end{frontmatter}

We dedicate this paper to the memory of Paul Joyce, friend and collaborator.

\section{Introduction}
In this paper we study aspects of the ancestral history of a random sample of $n$ DNA sequences, conditional on features of the haplotype configuration obtained at the present time, labeled 0. Initially, we assume an infinitely-many-sites mutation model with constant population size where the sites are in a completely linked region of DNA where there is no recombination. We begin with some theory that describes the effects of mutations between time 0 and time $t$ in the past. We describe the distribution of the quantities  $\widetilde{S}_n(t), A_n(t), \widetilde{K}_n(t)$ and $A_n^\theta(t)$, where $\widetilde{S}_n(t)$ is the number of mutations that have arisen in $(0,t)$, $A_n(t)$ is the number of ancestors at time $t$, $\widetilde{K}_n(t)$ is the number of distinct haplotypes  present in those ancestors that are formed by the mutations in $(0,t)$, and $A_n^\theta(t)$ is the number of  ancestors at time  $t$ whose descendants have no further mutation in $(0,t)$. 

Results such as these describe the effects of \emph{new} variation, that arising from the present time 0 to time $t$ in the past. We then describe the effects of standing variation by providing a simulation approach for studying the joint distribution of $S_n$ and $K_n$, the number of segregating sites and haplotypes in the sample at time 0 resp., and $S_n(t), K_n(t)$, the number of segregating sites and haplotypes, resp.,  in the $A_n(t)$ ancestors formed by standing variation arising after time $t$ in the past. This approach allows for variable population size, as well as essentially arbitrary binary branching models such as the Yule process.

If the full haplotype configuration and the number of segregating sites is available, then there are many other questions that may be asked about the ancestral history of the current sample. A sequential importance sampling algorithm for studying the ancestral history of a sample of genes conditional on  the haplotype configuration and the number of segregating sites is developed.
If the complete pattern of mutations on haplotypes were known, then a perfect phylogeny, a genetree, could be constructed and ancestral inference, such as the ages of mutations and time to the most recent common ancestor made conditional on the genetree topology. An example appears in the analysis in \cite{GT1994b}. 

Here we consider the case when just the haplotype frequencies and the number of mutations are known. The theory and  computational algorithms are then much simpler. The implementation computes, {\em inter alia},  the  probability of a sample configuration of haplotypes and segregating sites in a stationary population (this is an extension of the Ewens Sampling Formula), the average coalescence times, mutation times, allele loss times and allele ages back in time conditional on the current haplotype configuration and number of segregating sites, and the conditional average allele configuration and distribution of ancestors at time $t$ in the past. 
In the infinitely-many-sites model the haplotype configuration follows the infinitely-many-alleles model so the haplotype configuration in a sample is the same as the allele configuration and has a probability distribution of the Ewens Sampling Formula. 

\section{Ancestral distributions in the coalescent: theory}
We begin by setting some notation. We let $A_n(t)$ denote the number of ancestors of a sample of size $n$ time $t$ ago, and $A^\theta_n(t)$ the number of those ancestors whose descendants in $(0,t)$ have no further mutation. We let
$\widetilde{S}_n(t)$ be the number of segregating sites that arise as mutations in $(0,t)$, $\widetilde{K}_n(t)$ be the number of haplotypes in these ancestors that arise from mutations in $(0,t)$, and $S_n(t)$ the number of mutations in the sample arising after  time $t$ .

Coalescence times between events are denoted by $T_2,T_3,\ldots$ In the constant size coalescent model these are independent exponential random variables with rates ${2\choose 2}, {3\choose 2}, \ldots$~\citep{K1982}.  When considering non-mutant lineages back in time it is appropriate to have lineages lost by either mutation or coalescence. Times between these events are denoted by
$T_2^\theta,T_3^\theta,\ldots$. In the constant population size case these are independent exponential random variables with rates ${2\choose 2}+2\theta,{3\choose 2} +3\theta,\ldots$, $\theta$ being the scaled mutation rate. We denote by $f^\theta_{nk}(t)$  the density of $T_n^\theta+\cdots + T_k^\theta$. 

The reader is referred to \citet{E1972} for the distribution of haplotypes in a sample, the Ewens Sampling Formula; \citet{G1980}, eq. (9), for an original derivation of the distribution of $A_\infty^\theta(t)$ from looking back in time in a diffusion process; \citet{T1984} for an original derivation of $A_n^\theta(t)$, connections with the Kingman coalescent and review of ancestral lineage distributions; \citet{T1984,T2004} for an introduction to ancestral lineage distributions and ancestral inference in population genetics; and \citet{G2006} for a review and new representations for ancestral distributions.

The distribution of $A_n^\theta(t)$ is given by 
\begin{eqnarray}
\bbbp(A_n^\theta(t) = k) &:=& q_{nk}^{\theta}(t)
\nonumber \\
&=&\sum_{j=k}^n\rho_j^{\theta}(t)
(-1)^{j-k}
\frac{ (2j+ \theta - 1)(k+\theta)_{(j-1)}}{k!(j-k)!}\cdot\frac{n_{[j]}}{(n+\theta)_{(j)}},
\nonumber \\
\label{trnQ}
\end{eqnarray}
where $\rho_j^\theta(t) = e^{-\frac{1}{2}j(j+\theta-1)t}$; see \citep{T1984}, eq. (5.2). $\bbbp(A_n(t) = k)$ is given by setting $\theta=0$.  The formula also holds for $n=\infty$, where the interpretation is that of the whole infinite population coalescent.

The falling factorial moments of $A_n^\theta(t)$ are
\begin{equation}
\mathbb{E}\big [A_n^\theta(t)_{[r]}\big ] =
\sum_{k=r}^\infty\rho_k^\theta(t)(2k+\theta-1){k-1\choose r-1}(\theta+k)_{(r-1)}\frac{n_{[k]}}{(n+\theta)_{(k)}};
\end{equation}
see \citep{T1984}, p.13.
A simple rate argument establishes the correspondence between the distribution of coalescent times and the number of ancestors, namely
\begin{equation}
f_{nl}^{\theta}(t) = \frac{2}{l(l+\theta-1)}\bbbp(A_n^{\theta}(t) = l).
\label{fA:0}
\end{equation}
In a stationary population the probability generating function (pgf) of the number of segregating sites in a sample of $n$ genes is
\begin{eqnarray}
H_n(z) &=& \prod_{j=1}^{n-1}\Big (1 + \frac{\theta}{j}(1-z)\Big )^{-1}
\nonumber \\
&=& (n-1)\cdot \frac{\Gamma(n-1)\Gamma(\theta(1-z)+1)}{\Gamma(n+\theta(1-z))}
\nonumber \\
&=& \int_0^1x^{\theta(1-z)}(n-1)(1-x)^{n-2}dx.
\label{spgf:0}
\end{eqnarray}
Therefore in a stationary population,
\begin{equation}
\bbbp(S_n = k) = (n-1)\int_0^1(1-x)^{n-2}e^{-(-\theta\log x)}\frac{(-\theta \log x)^k}{k!}dx,
\label{sn:20}
\end{equation}
which is a Poisson mixture with mean $-\theta \log X$, where $X$ has density
$(n-1)(1-x)^{n-2}$, $0 < x < 1$. 
A calculation in \citet{T2004} (5.3.6) using the first line of (\ref{spgf:0}) shows that
\[
\bbbp(S_n = k)
= \frac{n-1}{\theta}\sum_{l=1}^{n-1}(-1)^{l-1}{n-2\choose l-1}\Bigg (\frac{\theta}{l+\theta}\Bigg )^{k+\theta}.
\]
 There is the simple recursion that
\begin{equation*}
n(n+\theta-1)\bbbp(S_n=s) = \theta \bbbp(S_n=s-1) + n(n-1)\bbbp(S_{n-1}=s),\> n \geq 2,
\end{equation*}
with $\bbbp(S_1=s) = \delta_{s0}$.

The distribution of the number of alleles in a sample of $n$ in a stationary population is
\begin{equation*}
\bbbp(K_n=k)= \theta^k |S_k^n|\,/\theta_{(n)},\> 1 \leq k \leq n,
\end{equation*}
where $\{S_k^n\}$ are Stirling numbers of the first kind.

\subsection{The joint distribution of $\widetilde{S}_n(t)$ and $A_n(t)$}
Measuring time back from the present, $(\widetilde{S}_n(t),A_n(t))$ is a Markov process beginning at $(0,n)$ at $t=0$  
such that
\begin{equation}
(s,l) \to \begin{cases}
(s+1,l)&\text{at~rate~}\theta l/2\\
(s,l-1)&\text{at~rate~} l(l-1)/2
\end{cases}.
\label{rates}
\end{equation} 
The pgf of the joint distribution of $A_n(t),\widetilde{S}_m(t)$ distribution is derived in \citet{G1981}, eqns (8), (9). The approach taken in the present paper is new.

The joint pgf of the number of mutations arising while there are $n,\ldots,k$ ancestors and the density of $T_n,T_{n-1},\ldots,T_{k}$ in an $n$-coalescent tree is
\begin{eqnarray}
&&\prod_{j=k}^ne^{\theta(z-1)jt_j/2}{j\choose 2}e^{-{j\choose 2}t_j/2}
\nonumber \\
&&=\frac{\prod_{j=k}^n{j\choose 2}}{\prod_{j=k}^n\frac{1}{2}j(j+\theta(1-z))}
\cdot
\prod_{j=k}^n\frac{1}{2}j(j+\theta(1-z)) e^{-\frac{1}{2}j(j-1+\theta(1-z))t_j}.
\end{eqnarray}
Integrating over $t_n+t_{n-1}+\cdots + t_k \leq t$, the pgf of the number of mutations on 
edges of the coalescent tree during times $T_n,\ldots,T_k$ and the distribution function of $T_n+\cdots +T_k$ is
\begin{eqnarray}
&&\frac{\prod_{j=k}^n{j\choose 2}}{\prod_{j=k}^n\frac{1}{2}j(j+\theta(1-z))}
\cdot \bbbp(T_k^{\theta(1-z)}+\cdots+T_n^{\theta(1-z)} \leq t)
\nonumber \\
&& =
\frac{\prod_{j=k}^n{j\choose 2}}{\prod_{j=k}^n\frac{1}{2}j(j+\theta(1-z))}
\bbbp(A_n^{\theta(1-z)}(t) \leq k)
\nonumber \\
&&= \prod_{j=k}^{n-1}\Big (1 + \frac{\theta(1-z)}{j}\Big )^{-1}\bbbp(A_n^{\theta(1-z)}(t) \leq k).
\end{eqnarray}
The joint probability that $A_n(t) = l$ and pgf of the number of mutations $\widetilde{S}_n(t)$  arising in $(0,t)$, which will be denoted by $G_l(z;t)$, is the probability that $T_n+\cdots +T_{l+1} = \tau < t$ (a coalescence necessarily occurs at $\tau$), there is no coalescence in $(\tau,t)$ and the pgf of the numbers of mutations on $T_n,\ldots, T_{l+1}$ and the $l$ lines from $\tau$ to $t$.
Therefore, with notation $\theta_z = \theta(1-z)$,
\begin{eqnarray}
G_l(z;t) &=& \frac{\prod_{j=l+1}^n{j\choose 2}}{\prod_{j=l+1}^n\frac{1}{2}j(j-1+\theta(1-z))}
\cdot \int_0^t e^{-{l\choose 2}(t-\tau) + \theta (z-1)l(t-\tau)}
f^{\theta(1-z)}_{n,l+1}(\tau)d\tau
\nonumber \\
&=&
\frac{\prod_{j=l+1}^n{j\choose 2}}{\prod_{j=l}^n\frac{1}{2}j(j-1+\theta(1-z))}
\cdot
f^{\theta(1-z)}_{n,l}(t)
\nonumber \\
&=&
\frac{\prod_{j=l+1}^n{j\choose 2}}{\prod_{j=l+1}^n\frac{1}{2}j(j-1+\theta_z)}
\cdot
\bbbp(A_n^{\theta_z}(t) = l)
\label{id:0}\\
&=&
 \prod_{j=l}^{n-1}\Big (1 + \frac{\theta_z}{j}\Big )^{-1}
\bbbp(A_n^{\theta_z}(t) = l)
\nonumber \\
&=&
 \prod_{j=l}^{n-1}\Big (1 + \frac{\theta_z}{j}\Big )^{-1}
\nonumber \\
&&\times
\sum_{j=l}^n\rho_j^{\theta_z}(t)
(-1)^{j-l}
\frac{ (2j+ \theta_z - 1)(l+\theta_z)_{(j-1)}}{l!(j-l)!}\cdot\frac{n_{[j]}}{(n+\theta_z)_{(j)}}.
\nonumber \\
\label{pgfprob:0}
\end{eqnarray} 
The identity (\ref{fA:0}) is used in the calculation.
As a check of (\ref{pgfprob:0}), if $z=1$ the probability that there are $l$ ancestors at $t$ is $\bbbp(A_n^0(t) = l)$, as it should be.
\ignore{
\st{\textit{Bob: I do not understand the next bit ... indices adrift??}  The marginal pgf of the number of mutations arising while $l$ to $j$ ancestors found by summing over $l$ in (\ref{pgfprob:0}) is 
\[
\prod_{j=l}^{n-1}\Big (1 + \frac{\theta_z}{j}\Big )^{-1}.
\] 
}
This agrees with a calculation similar to (\ref{sn:20}) where 
\[
\prod_{j=l}^{n-1}\Big (1 + \frac{\theta_z}{j}\Big )^{-1}
= \frac{\Gamma (n)}{\Gamma(l)\Gamma (n-l)}\cdot \frac{\Gamma(l+\theta_z)\Gamma(n-l)}{\Gamma(n+\theta_z)}
\]
is the pgf of a Poisson mixture with a rate $-\theta\log X$, where $X$ has a Beta$(l,n-l)$ distribution. That is, the pgf is
\begin{equation}
B(l,n-l)^{-1}
\int_0^1x^{l-1}(1-x)^{n-l-1}e^{-(-\theta\log x)}\frac{(-\theta \log x)^k}{k!}dx,
\label{nldist:0}
\end{equation}
for $k=0,1,\ldots$. 
The pgf of the total number of mutations accumulating while there are greater than or equal to $l$ ancestors is
\begin{equation}
\frac{l(l-1)}{2}\int_0^\infty G_l(z;t)dt = \frac{l-1}{(l-1+\theta_z)}
\cdot \prod_{j=l}^{n-1}\Big (1 + \frac{\theta_z}{j}\Big )^{-1}
=\prod_{j=l-1}^{n-1}\Big (1 + \frac{\theta_z}{j}\Big )^{-1}
\label{total:0}
\end{equation}
}

The marginal pgf of the number of mutations arising while there are at least $l$ ancestors  is given by
\begin{equation}
\prod_{j=l-1}^{n-1}\Big (1 + \frac{\theta_z}{j}\Big )^{-1}.
\label{pgf:1000}
\end{equation} 
To see this, calculate the pgf of the total number of mutations accumulating while there are at least $l$ ancestors as
\begin{equation}
\frac{l(l-1)}{2}\int_0^\infty G_l(z;t)dt = \frac{l-1}{(l-1+\theta_z)}
\cdot \prod_{j=l}^{n-1}\Big (1 + \frac{\theta_z}{j}\Big )^{-1}
=\prod_{j=l-1}^{n-1}\Big (1 + \frac{\theta_z}{j}\Big )^{-1}.
\label{total:0}
\end{equation}
In a calculation similar to (\ref{sn:20}), 
\[
\prod_{j=l-1}^{n-1}\Big (1 + \frac{\theta_z}{j}\Big )^{-1}
= \frac{\Gamma (n)}{\Gamma(l-1)\Gamma (n-l+1)}\cdot \frac{\Gamma(l-1+\theta_z)\Gamma(n-l+1)}{\Gamma(n+\theta_z)}
\]
is the pgf of a Poisson mixture with a rate $-\theta\log X$, where $X$ has a Beta$(l-1,n-l+1)$ distribution. Hence the probability of $k$ mutations, from the pgf (\ref{pgf:1000}), is
\begin{equation}
B(l-1,n-l+1)^{-1}
\int_0^1x^{l-2}(1-x)^{n-l}e^{-(-\theta\log x)}\frac{(-\theta \log x)^k}{k!}dx,
\label{nldist:0}
\end{equation}
for $k=0,1,\ldots$.

The pgf of the number of mutations arising in $(0,t)$ (not counting mutations when there is one ancestor) is
$
\sum_{l=2}^n G_l(z;t).
$
The joint pgf for $S_n(t)$, $\widetilde{S}_n(t)$ and the probability of $l$ ancestors time $t$ ago is $Q_l(r)G_l(z;t)$, where $Q_l(r)$ is the pgf of the number of segregating sites in a sample of $l$ from the population at the initial time. In a stationary population
$Q_l(r) = H_l(r)$, defined in (\ref{spgf:0}).

The probability $\bbbp(\widetilde{S}_n(t) = k)$ does not have a simple form, but if one considers standing variation in a stationary population as well as mutations in $(0,t)$ then the distribution of the number of mutations is as in a stationary population at time $t$ and the conditional distribution of the number of ancestors given the number of mutations is easier. The pgf/probability of the number of mutations and number of ancestors at time $t$ in the past is then
\begin{equation}
 G^*_l(z,t) = \prod_{j=2}^{n-1}\Big (1 + \frac{\theta_z}{j}\Big )^{-1}
\bbbp(A_n^{\theta_z}(t) = l).
\label{allmutations:0}
\end{equation}
Summing (\ref{allmutations:0}) over $l=2,\ldots,n$, we see that the marginal pgf of the number of mutations  is
\begin{equation}
\prod_{j=2}^{n-1}\Big (1 + \frac{\theta_z}{j}\Big )^{-1},
\label{check:0}
\end{equation}
and the distribution of the number of ancestors at time $t$ back, found by setting $z=1$ in (\ref{allmutations:0}), is
\begin{equation}
\bbbp(A_n^{0}(t) = l).
\label{check:1}
\end{equation}
Equations (\ref{check:0}) and (\ref{check:1}) recover earlier results.
The joint pgf and expected number of ancestors time $t$ ago is
\begin{eqnarray}
&& \prod_{j=2}^{n-1}\Big (1 + \frac{\theta_z}{j}\Big )^{-1}
 \mathbb{E}\big [A_n^{\theta_z}(t) \big ]
 \nonumber \\
 &&~= 
 \prod_{j=2}^{n-1}\Big (1 + \frac{\theta_z}{j}\Big )^{-1}
 \sum_{k=1}^n\rho_k^{\theta_z}(t)(2k+\theta_z-1)\frac{n_{[k]}}{(n+\theta_z)_{(k)}}.
 \label{allmutations:1}
\end{eqnarray}
Inversion of (\ref{allmutations:1}) is straightforward, if a little messy.
We calculate $\mathbb{E}\big [A_n(t)\mid S_n = r\big ]$.
 Let $a(r,k)$ be the coefficient of $z^r$ in
\[
M(k,z)= \prod_{j=2}^{n-1}\Big (1 + \frac{\theta_z}{j}\Big )^{-1}
 (n+\theta_z)_{(k)}^{-1}.
 \]
 Then by equating coefficients of $z^r$ in 
\[
(n+\theta+k-1-\theta z)M(k,z) = M(k-1,z),
\]
 for $k=1,2,\ldots, n$, $r = 0,1,\ldots$, we get
 \[
 (n+\theta+k-1)a(r,k) = \theta a(r-1,k) + a(r,k-1).
 \]
 Note that $a(r,0) = \bbbp(S_n=r)$.
 Let 
 \[
 b(r,k) = (2k+\theta-1)a(r,k) - \theta a(r,k-1)
 \]
and finally
 \[
 c(r,k) = \sum_{m=0}^r e^{-k\theta t/2}\frac{(k\theta t/2)^m}{m!}b(r-m,k).
 \]
 Then 
\begin{equation}
\mathbb{E}\big [A_n(t)\mid S_n = r\big ] 
=  \frac{\sum_{k=1}^n\rho_k(t)c(r,k) n_{[k]}}{\bbbp\big (S_n = r\big )}.
\label{allmutations:3}
\end{equation}
The pgfs (\ref{pgfprob:0}) and (\ref{allmutations:0}) and the formula (\ref{allmutations:3}) are new.

In a later section we give a rejection algorithm for simulating the distribution of $A_n(t)$ given $S_n=r$.


\subsection{The joint distribution of $\widetilde{S}_n(t)$, $\widetilde{K}_n(t)$, $A_n(t),A_n^\theta(t)$}
A method for computing the stationary joint distribution of $(S_n,K_n)$ is derived in \cite{G1982}, exploiting  a diffusion generator in a model with $K$ alleles, then letting $K\to \infty$.  It is found numerically that the joint distribution is strongly diagonal in that, approximately, $S_n=K_n-1$.  It is also shown that  $S_n - K_n + 1$ has a proper limit distribution when $n \to \infty$, even though both $S_n$ and $K_n$ tend to infinity. It is always true that $\widetilde{S}_n(t) - \widetilde{K}_n(t) + 1 \geq 0$. 

Here we consider the joint distribution of $(\widetilde{S}_n(t),\widetilde{K}_n(t),A_n(t),A_n^\theta(t))$ using a coalescent treatment different  from that of \citet{G1982}. To obtain a Markov process, consider  $(\widetilde{S}_n(t), \widetilde{K}_n(t), A_n(t), A_n^\theta(t)), t \geq 0$, beginning from $(0,0,n,n)$. Let $B_n(t) = A_n(t) - A_n^\theta(t)$. Then 
\begin{equation}
(s,k,b,a^\theta) \to \begin{cases}
(s+1,k,b,a^\theta)&\text{at~rate~}\theta b/2\\
(s+1,k+1,b+1,a^\theta-1)&\text{at~rate~} \theta a^\theta/2\\
(s,k,b-1,a^\theta)&\text{at~rate~} \big (b(b-1) + 2 ba^\theta\big )/2\\
(s,k,b,a^\theta-1)&\text{at~rate~} a^\theta(a^\theta-1)/2
\end{cases}
\label{rates:ska:0}
\end{equation}
The total coalescence rate is $a(a-1)/2=(a^\theta+b)(a^\theta+b-1)/2$. The total mutation rate is $a\theta/2=\theta(a^\theta+b)/2$. This process counts mutations and alleles as they arrive back in time from $t$ in a sample of $n$ and its ancestors, by considering the two groups of lines $a-a^\theta,a^\theta$ and in which groups mutations or coalescences occur.

The simpler process $(A_n(t),A_n^\theta(t))$ has rates
\begin{equation}
(a,a^\theta) \to \begin{cases}
(a-1,a^\theta)&\text{at~rate~} (b(b-1)+2ba^\theta)/2\\
(a-1,a^\theta-1)&\text{at~rate~} a^\theta(a^\theta-1)/2\\
(a,a^\theta-1)&\text{at~rate~} a^\theta \theta/2\\
\end{cases}
\label{rates:ska:10}
\end{equation}
with total rate of $\left(a(a-1)+a^\theta\right)/2$.
The marginal transition functions are known explicitly from (\ref{trnQ}) when $\theta=0$ and $\theta>0$.

A slightly different approach, described in equation (2.9) of \citet{G1981}, is to consider sample paths that have $i$ mutations and $j$ alleles at $t$ starting from fixed $m=a^\theta$, $r=a$ as being stationary distributions which satisfy the following recursive system:
\begin{eqnarray}
&&a(a-1+\theta)p(i,j;a^\theta,a)
\nonumber \\
&&~~=(a-a^\theta)\theta p(i-1,j;a^\theta,a) + (a+a^\theta-1)(a-a^\theta)p(i,j;a^\theta,a-1)
\nonumber \\
&&~~~~+a^\theta\theta p(i-1,j-1;a^\theta-1,a) + a^\theta(a^\theta-1)p(i,j;a^\theta-1,a-1),
\label{Gks:0},
\end{eqnarray}
for $i=0,1,\ldots,a^\theta;j=1,2\ldots\text{~and~}a \geq a^\theta \geq 2.$ The boundary probabilities satisfy
\begin{eqnarray}
p(i,j;1,a) &=& 0,\>j>1
\nonumber \\
(a+\theta-1)p(i,1;1,a) &=& (a-1)p(i,1;1,a-1) + \theta p(i,1;1,a),\>a=2,3,\ldots
\nonumber \\
p(i,1;1;1) &=& \delta_{i0}.
\label{Gks:1}
\end{eqnarray}
Then we want $p(i,j;n,n)$, which can be computed from (\ref{Gks:0}). These equations can be argued from a probabilistic perspective, though a different approach is taken in \citet{G1981}, where the identification $m=a^\theta$, $r=a$ is not made.
(\ref{Gks:0}) is analogous, but more complex, to a similar recursion for the number of mutations which could be written as
\begin{equation}
a(a+\theta-1)p(i;a) = \theta ap(i-1;a) + (a-1)p(i;a-1).
\end{equation}
for $i=0,1 \ldots,\>a=2,3,\ldots$.

A pgf version of (\ref{Gks:0}) for
\[
q(a^\theta,a) = \sum_{i=0}^\infty\sum_{j=1}^{i+1}p(i,j;a^\theta,a)z^iw^j
\]
appears in   \cite{G1982}, eq. (2.8).

\ignore{
\textcolor{blue}{We may or may not want this next paragraph in the final version}\\
In a time interval of length $\tau$, if the initial configuration is $(s,k,b,a^\theta)$,
then the probability that the configuration is $(s_\tau,k_\tau,b_\tau,a^\theta_\tau)$  afer time $\tau$, when there are no coalescences, is now calculated.
Suppose that there are $\alpha$ mutations on mutant lineages $b$ and $\beta$ mutations on non-mutant lineages $a^\theta$, $\alpha + \beta = s_\tau - s$. The probability that there are $a^\theta_\tau$ non-mutant lineages remaining after $\beta$ mutations in the group is given by inclusion-exclusion. The probability that a specific $r$ lineages are non-mutant is
\[
e^{-\theta r\tau/2}\big ( 1 - e^{-\theta(a^\theta-r)\tau/2}\big )^\beta,
\]
so the probability that exactly $a_\tau^\theta$ of $a^\theta$ lineages are non-mutant is
\begin{equation}
p(a_\tau^\theta\mid a^\theta, \beta) =\sum_{r=0}^{a^\theta}{a^\theta\choose r}{r\choose a_\tau^\theta}(-1)^{r-a_\tau^\theta}
e^{-\theta r\tau/2}\big ( 1 - e^{-\theta(a^\theta-r)\tau/2}\big )^\beta.
\end{equation}
Note that $b_\tau = b + a - a_\tau$, $k_\tau = k + a-a_\tau$. Then
\begin{equation}
p((s_\tau,k_\tau,b_\tau,a^\theta_\tau)\mid (s,k,b,a^\theta)) = \sum_{\beta}\frac{\big (\theta b\tau/2\big )^{s_\tau-s-\beta}}
{(s_\tau-s-\beta)!}
e^{-\theta b\tau/2}
p(a_\tau^\theta\mid a^\theta, \beta).
\end{equation}
}

\section{Simulation-based approaches}

In this section we develop theory that will be used in the rejection sampling and importance sampling schemes. We begin by recalling a simulation method that generates stationary samples of haplotype counts together with the number of mutations in the tree. This works for constant population size coalescent models.

\subsection{Growing a tree} A useful way to simulate the ancestral history of haplotype configurations in age order together with the mutations  is to use a condensation of an algorithm for \emph{growing} a gene tree whose nodes are mutations. The algorithm for the tree, described in~\citet{EG1987},Theorem 5.4 and \citet{G1989}, p.7, is the following.
\bigskip

 \emph{Algorithm 1}
\begin{enumerate}
\item
Start with a tree of two leaves (individuals) as two edges joined at a node.
\item
When there are $m$ leaves, select one of the $m$  leaves at random and duplicate from the same immediate mutation node with probability $(m-1)/(\theta+m-1)$;
or add a mutation node on the chosen edge with probability $\theta/(\theta+m-1)$.
\item
To get a sample of $n$, stop when there are first $n+1$ leaves and select the configuration just before the last leaf appeared.
\end{enumerate}
Times between events can be added as exponential random variables
$T^\theta_2,\ldots, T^\theta_n$.

If we just look at haplotype frequencies in age order from the oldest and keep track of $s$, the number of accumulated mutations,  the state space is $(m_1,\ldots,m_k;s)$ and transitions are Markovian. The condensed algorithm is the following.
\bigskip

\emph{Algorithm 2}
\begin{enumerate}
\item Start with a configuration $(n_1=2;0)$ of two identical oldest haplotypes and no mutations. 
\item
When there are $m$ individuals make a transition $(m_1,\ldots,m_k;s) \to (m_1,\ldots,m_j+1,\ldots ,m_k;s)$ with probability $(m_j/m)(m-1)/(m+\theta-1)$; or 
$(m_1,\ldots,m_k;s) \to (m_1,\ldots, m_l-1,\ldots,m_k,1;s+1)$ with probability
$(m_l/m)\theta/(m+\theta-1)$.
\item
To get a sample of $n$, stop when there are first $n+1$ leaves and select the configuration just before the last leaf appeared.
\end{enumerate}
This  algorithm is useful for simulation of an ancestral path forward in time which contains full information of haplotype count configurations and mutations.

Let $p^\circ(\bm{m})/\prod_{j=1}^m\alpha_j!$ be the probability of a non-age labelled configuration $\bm{m}$ (labelled in an arbitrary order), with $\alpha_j$ the number of allele frequencies equal to $j$, under the Markov chain without the stopping rule.
Then
\begin{eqnarray}
p^\circ(\bm{m};s) &= &
\frac{m-2}{m+\theta-2}\sum_{m_j>1}\frac{m_j-1}{m-1}p^\circ(m_1,\ldots,m_j-1,\ldots,m_k;s)
\nonumber \\
&&+\frac{\theta}{m+\theta-1}\sum_{i,l:m_i=1}\frac{m_l+1-\delta_{li}}{m}
p^\circ(\bm{m}+\bm{e}_l-\bm{e}_i;s-1).
\label{pcirc:0}
\end{eqnarray}
The sample probability $p(\bm{n};s)/\prod_{j=1}^n\alpha_j!$ is such that
\[
p(\bm{n};s) = \frac{n-1}{n+\theta-1}p^\circ(\bm{n};s).
\]
A recursion is therefore
\begin{eqnarray}
p(\bm{n};s) &= &
\frac{n-1}{n+\theta-1}\sum_{n_j>1}\frac{n_j-1}{n-1}p(\bm{n}-\bm{e}_j;s)
\nonumber \\
&&+\frac{\theta}{n+\theta-1}\sum_{i,l:n_i=1}
\frac{n_l+1-\delta_{li}}{n}
p(\bm{n}-\bm{e}_i+\bm{e}_l;s-1).
\label{pcirc:1}
\end{eqnarray}
The recursion in (\ref{pcirc:1}) is used in calculating the importance sampling weights in Section 4.

Note that the last sum includes the case when $i=l$, and it can be written as
\[
\sum_{i\ne l:n_i=1}
\frac{n_l+1}{n}
p(\bm{n}-\bm{e}_i-\bm{e}_l;s-1)
+ \alpha_1\frac{1}{n}p(\bm{n};s-1).
\]
If $a,s$ are the number of haplotypes and mutations, respectively, with $\bm{n}$ then $s -a +1 \geq 0$ so $p(\bm{n};s) = 0$ if $s -a +1 <0$ in the recursive equations.
The Ewens Sampling formula~\citep{E1972} satisfies a similar equation to (\ref{pcirc:1}), summing over $s$, with
\[
p(\bm{n};\cdot) = \frac{n!}{n_1\cdots n_k}\cdot \frac{\theta^k}{\theta_{(n)}}.
\]
The probability of the sample configuration is
\[
p(\cdot)/\prod_{j=1}^n\alpha_j! = \frac{n!}{\prod_{j=1}^kn_j\alpha_j!}\cdot \frac{\theta^k}{\theta_{(n)}} = \frac{n!}{\prod_{j=1}^kj^{\alpha_j}\alpha_j!}\cdot \frac{\theta^k}{\theta_{(n)}}.
\] 

The joint distribution of $(\bm{n};s)$ is studied in \citet{IZMPTR2005}, where their derived equation (1) is essentially the same as (\ref{pcirc:1}). The authors then use these recursive equations to derive exact probabilities for small sample sizes. In Section \ref{S:4} we use an importance sampling technique based on (\ref{pcirc:1}) which allows calculation for much larger sample sizes than are possible with exact calculation and we also carry out ancestral inference back in time.

\subsection{The number of ancestors and mutations}

The simplest ancestral question is to ask about the distribution of $A_n(t)$ conditional on $S_n=s$ segregating sites in the sample of size $n$. Rejection algorithms work well for  problems like this, as was illustrated by \citet{TBGD1997}.
They took a general Bayesian approach in which $\theta$ was considered as a random variable and the times $T_n,\ldots,T_2$ come from a variable population size coalescent model.  The combinatorics for a general binary coalescent tree are studied in \cite{GT1998}, where the sample frequency spectrum and the mean age of a mutation known to have $b$ descendents in a sample of $n$ genes are calculated.

We illustrate by describing how to generate observations from the conditional distribution of $(\theta, A_n(t))$ conditional on $S_n = s$.  To this end, let $W_1=T_n, W_2=T_n+T_{n-1},\ldots, W_j = T_n+\cdots +  T_{n-j+1},\ldots,
W_{n-1} = T_n + \cdot + T_2$. 
$W_{n-1}$ is the height of the coalescent tree. Define the total length of the tree as $L_n = nT_n+\cdots+2T_2$. Form a set of bins as follows:
\[
B_1=(0,W_1],\ldots,B_{n-1}=(W_{n-2},W_{n-1}],B_n=(W_{n-1},\infty).
\]
Let $J$ be the bin that covers $t$. Then the number of ancestors at time $t$ is
\begin{equation}
A_n(t) = n-J+1,
\label{anc_bin:0}
\end{equation}
and the length of the coalescent tree from 0 back to $t$ is
\begin{equation}
\widetilde{L}_n(t) = \begin{cases}
nt&\text{if~}J=1\\
L_n&\text{if~}J=n\\
\sum_{l=n-J+2}^nlT_l + (n-J+1)(t-W_{J-1})&\text{if~}2 \leq J \leq n-1.
\end{cases}
\label{length:1}
\end{equation}
Finally, define the length of the coalescent tree from time $t$ to the most recent common ancestor as $L_n(t) = L_n-\widetilde{L}_n(t)$. Let $\widetilde{S}_n(t)$ be the number of mutations arising in $(0,t)$ and let $S_n(t)$ the number of mutations from $t$ to the most recent common ancestor. Then, conditional on $T_n,\ldots,T_2$ and $J$,
\[
\widetilde{S}_n(t) \sim {\rm Po}\Big (\theta\widetilde{L}_n(t)/2\Big ),\> S_n(t) \sim {\rm Po}\Big (\theta L_n(t)/2\Big),
\]
where ${\rm Po}(\lambda)$ denotes the Poisson distribution with parameter $\lambda$,
such that
\[
Po(\lambda)\{s\} = \frac{e^{-\lambda}\lambda^s}{s!},\>s = 0,1,\ldots,
\]
and $\widetilde{S}_n(t)$ and $S_n(t)$ are conditionally independent. The total number of segregating sites in the sample at time 0 is $S_n = \widetilde{S}_n(t) + S_n(t)$. The simplest algorithm gives the probability distribution of $J$, and therefore the distribution of $A_n(t)$, conditional on $S_n=s$. 
\bigskip

\emph{Algorithm 3}\label{Algorithm:3}
\begin{enumerate}
\item
Simulate $\theta$ from the prior, $\pi(\cdot)$
\item
Simulate $T_n,\ldots,T_2$ from an appropriate coalescent model
\item
Compute $J$, $A = A_n(t) = n - J + 1$, $L_n$
\item
Accept $(\theta,A)$ as an observation from the posterior with probability
\[
h = \frac{{\rm Po}(\theta L_n/2)\{s\}}{{\rm Po}(s)\{s\}}
\]
\end{enumerate}

We also note that the same rejection approach may be used to approximate conditional distributions for many other ancestral variables. For example, to study the distribution of the number of mutations $S_n(t)$ present in the ancestors at time $t$, which is a measure of the standing variation at time $t$, we can use the following.
\bigskip

\emph{Algorithm 4}\label{Algorithm:4}
\begin{enumerate}
\item
Simulate $\theta$ from the prior, $\pi(\cdot)$.
\item
Simulate $T_n,\ldots,T_2$ from an appropriate coalescent model.
\item
Compute $J$, $A = A_n(t) = n - J + 1$, $L_n(t)$ and $\widetilde{L}_n(t)$, defined in (\ref{length:1}).
\item 
Accept $(\theta,A)$ with probability
\[
h = \frac{{\rm Po}(\theta L_n/2)\{s\}}{{\rm Po}(s)\{s\}}
\]
and else return to 1.
\item Simulate $S$ from a Binomial($s,L_n(t)/L_n$) distribution, and return $(\theta,A,S)$ as an observation from the posterior of $(\theta,A_n(t), S_n(t))$ given $S_n = s$.
\end{enumerate}

We may treat $\theta$ as fixed in this approach (that is, as having a degenerate prior), the approach then addressing the problems studied in the first section of the paper. 

\citet{BR2007} construct a rejection algorithm for maximum likelihood estimation of the number of ancestral lineages at time $t$ back based on the sample frequency spectrum. The algorithm is constructed differently from Algorithms 3 and 4.

\subsection{Hammer et al. example}\label{HammSect}
As an illustration we consider the Y  chromosome data of 1544 sequences from \citet{H1998}. In this paper a perfect phylogeny was constructed from the sequence data and the program GENETREE was used to find the TMRCA and ages of mutations in the ancestral tree shown in Figure 7 of \citep{H1998} with $\theta=2.5$. There were 9 segregating sites and 10 haplotypes observed in the data. The unconditional expected height of the coalescent tree is 2 time units, and the expected height conditional on $s = 9$ segregating sites is 1.21 units, this latter found from the method in~\citep{TBGD1997}.

We generated 10,000  repetitions of the previous algorithm for a series of times $t$, obtaining the information in Table~\ref{Table1} for the conditional expectations of $S_n(t)$ and $A_n(t)$.

\begin{table}[htb]\label{Table:1}
\begin{center}
\begin{tabular}{|c|cc|cc|}
\hline
$t$ & $S_n(t)$ &SE & $A_n(t)$ & SE\\
\hline
0 & 9.0 & & 1544 & \\
0.1 & 3.18 & 0.031 & 19.47 & 0.051 \\
0.5 & 1.11 & 0.023 & 3.72 & 0.023 \\
1.0 & 0.38 & 0.015 & 1.78 & 0.015 \\
1.5 & 0.12 & 0.009 & 1.25 & 0.009 \\
\hline
\end{tabular}
\caption{Result of 10,000 runs for the constant population size coalescent model with a sample of size $n = 1,544$ and $s = 9$ segregating sites. Table shows average value of $A_n(t), S_n(t)$ given $S_n = 9$. Righthand columns give SE of the mean.}\label{Table1}
\end{center}
\end{table}

The simple rejection schemes illustrated here are not as useful for considering more complicated summaries of the data. In the next section we show how to exploit an importance sampling approach to derive conditional distributions given the haplotype frequency distribution and the number of mutations.

\section{Importance sampling}\label{S:4}

Sequential importance sampling for ancestral inference in population genetics has a long history, illustrated by  \citet{GT1994a,GT1994b,GT1999,F1999,SD2000,G2002,DG2004}. The technique can be described as constructing a proposal distribution for events back in time, simulating back in time, then correcting for the approximate proposal distribution by calculating the exact probability of the path forward in time and taking the ratio of the probability of the forward path divided by the approximate probability of the backward path as the importance weight. If there are $r$ simulation runs then an empirical ancestral history is returned as 
$(\widehat{p}_1,{\cal H}_1), \ldots, (\widehat{p}_r,{\cal H}_r)$  where $\{\widehat{p}_j\}$ are the importance weights scaled to add to unity and $\{{\cal H}_j\}$ are the histories.
A general reference to sequential importance sampling is \citet{L2001}.

Choosing a proposal distribution is an art. We use the principal of choosing a lineage which can be involved in a transition back in time uniformly. This has a theoretical justification, described in  \citep{SD2000,DG2004}. Sequential importance sampling for a haplotype configuration can be regarded as a simplification of the technique used for a complete genetree, constructed as a perfect phylogeny from the pattern of mutations on DNA sequences. The simpler scheme counts different haplotypes of the sequences, the extra information being the number of mutations back to the most recent common ancestor. Time information, such as coalescence times, ages of mutations and time to the most recent ancestor can be included. Details of how to include time in an importance sampling algorithm are in \citet{GT1997}. We develop a new importance sampling approach for the Kingman coalescent models conditional on an observed configuration $\bm{n},s$.


The proposal distribution $\widehat{p}$ for reverse transitions in a haplotype history is detailed in the following equations. Suppose a current configuration is $\bm{n}=(n_1,\ldots,n_k)$, the number of mutations to the most recent common ancestor is $s$, and the number of singletons is $q$.

For $k>2$, if $s-k+1>0$,
\begin{eqnarray}
\widehat{p}(\bm{n}-\bm{e_i};s\mid \bm{n};s) &=& \frac{n_i}{n}\text{~if~} n_i > 1
\label{nosingleton}\\
\widehat{p}(\bm{n}-\bm{e}_i + \bm{e}_l;s-1\mid \bm{n};s) &=& \frac{n_l}{n}\cdot\frac{1}{n},
\>n_i=1, l\ne i
\nonumber \\
\widehat{p}(\bm{n};s-1\mid \bm{n};s) &=&\frac{q}{n}\cdot \frac{1}{n},
\label{singleton}
\end{eqnarray}
or if $s-k+1=0$ then
\begin{eqnarray}
\widehat{p}(\bm{n}-\bm{e_i};s\mid \bm{n};s) &=& \frac{n_i}{n}\text{~if~} n_i > 1
\nonumber\\
\widehat{p}(\bm{n}-\bm{e}_i + \bm{e}_l;s-1\mid \bm{n};s) &=& \frac{n_l}{n-1}\cdot\frac{1}{n},
\>n_i=1, l\ne i.
\label{second:0}
\end{eqnarray}
The first factors in (\ref{singleton}) involve a choice of either mutations that define allele types and those which appear on lineages between defined alleles. 
Importance weights for transitions back in time are therefore
\begin{eqnarray*}
(\bm{n},s) \to (\bm{n}-\bm{e}_i,s),n_j>1:&& 
\frac{n_i-1}{n+\theta-1}\cdot 
\frac{1}{\widehat{p}(\bm{n}-\bm{e}_i,s\mid \bm{n},s)}
\nonumber \\
(\bm{n},s) \to (\bm{n}-\bm{e}_i+\bm{e}_l,s-1),n_k=1,k\ne l:&&
 \frac{\theta}{n+\theta-1}\cdot \frac{n_l+1}{n}
\frac{1}{\widehat{p}(\bm{n}- \bm{e}_i+\bm{e}_l,s-1\mid \bm{n},s)}
\nonumber \\
(\bm{n},s) \to (\bm{n},s-1):&&
\frac{\theta}{n+\theta-1}\cdot\frac{1}{n} 
\frac{1}{\widehat{p}(\bm{n},s-1\mid \bm{n},s)}
\end{eqnarray*}
When $k=2$ we have to consider the following possible cases:\\
(a) If $n_1>1,n_2>1$
\[
\widehat{p}(\bm{n}-\bm{e}_i;s\mid \bm{n};s) = \frac{n_i}{n},\>i=1,2,
\] 
(b) if $n_1>1,n_2=1,s>1$, 
\begin{eqnarray*}
\widehat{p}(\bm{n}-\bm{e}_1;s\mid \bm{n};s) &=& \frac{n_1}{n},\\
\widehat{p}(\bm{n};s-1\mid\bm{n};s) &= &\frac{1}{n},
\end{eqnarray*}
(c) if $n_1>1,n_2=1,s=1$,
\begin{eqnarray*}
\widehat{p}(\bm{n}-\bm{e}_1;1\mid \bm{n};1) &=& \frac{n_1}{n},\\
\widehat{p}(\bm{n}+\bm{e}_1-\bm{e}_2,0\mid \bm{n};1) &=& \frac{1}{n},
\end{eqnarray*}
(d) similarly when $n_1=1,n_2>1$,\\
(e) if $n_1=1,n_2=1,s>1$,
\[
\widehat{p}(\bm{n};s-1\mid \bm{n};s) = 1,
\]
(f) if $n_1=1,n_2=1,s=1$
\[
\widehat{p}(\bm{n}+\bm{e}_1-\bm{e}_2;0\mid \bm{n};1) = 1.
\]
Importance weights are:
\begin{eqnarray*}
&&\text{(a)~} (\bm{n},s) \to (\bm{n}-\bm{e}_i,s):\frac{n_i-1}{n+\theta-1}\cdot 
\frac{1}{\widehat{p}(\bm{n}-\bm{e}_i,s\mid \bm{n},s)},\>n_1,n_2>1\\
&&\text{(b)~} (\bm{n},s) \to (\bm{n}-\bm{e}_1,s):\frac{n_1-1}{n+\theta-1}\cdot
\frac{1}{\widehat{p}(\bm{n}-\bm{e}_1,s\mid \bm{n},s)},\>n_1>1,n_2=1,s>1\\
&&\text{\phantom{(b)~}} (\bm{n},s) \to  (\bm{n},s-1): \frac{\theta}{n+\theta-1}\cdot\frac{1}{n}\cdot
\frac{1} {\widehat{p}(\bm{n},s-1\mid \bm{n},s)}\\
&&\text{(c)~}(\bm{n},1) \to (\bm{n}-\bm{e}_1,1):\frac{n_1-1}{n+\theta-1}\cdot
\frac{1}{\widehat{p}(\bm{n}-\bm{e}_1,1\mid \bm{n},1)},\>n_1>1,n_2=1,s=1\\
&&\text{\phantom{(c)~}}(\bm{n},1) \to (\bm{n}+\bm{e}_1-\bm{e}_2,0):
\frac{\theta}{n+\theta-1}\cdot\frac{1}{\widehat{p}(\bm{n}+\bm{e}_1-e_2,1\mid \bm{n},1)}\\
&&\text{(e)~}(\bm{n},s) \to  (\bm{n},s-1):
\frac{\theta}{1+\theta}\cdot\frac{1}{\widehat{p}(\bm{n},s-1\mid \bm{n},s)},\>n_1=1,n_2=1,s>1\\
&&\text{(f)~}(\bm{n},1) \to (\bm{n}+\bm{e}_1-\bm{e}_2,0):\frac{2\theta}{1+\theta}\cdot\frac{1}{\widehat{p}(\bm{n}+\bm{e}_1-e_2,0\mid \bm{n},1)},\>n_1=1,n_2=1,s=1.
\end{eqnarray*}

\subsection{Implementation}
Our implementation provides
\begin{itemize}
\item
The probability of a sample configuration of haplotypes and number of segregating sites. This is an extension of the Ewens Sampling Formula, which is the probability of the configuration of haplotypes.
\end{itemize}
The next calculations are conditional on the configuration of haplotypes and segregating sites at  time 0.
\begin{itemize} 
\item
The average coalescence times in the past. 
\item
The average mutation times in the past. 
\item
The average times when alleles are lost in the past. (The time of loss of the last haplotype  is truncated at the TMRCA if not lost by mutation.)
\item
The average allele ages in the past. 
\item
The distribution of ancestor lines and the average allele configuration at a given time  in the past.
\end{itemize}

The program also implements a variable population size option with exponential growth. Coalescent times then have a distribution that depends on the time when they occur. We do not go into detail here, but refer the reader to  \citet{GT1994c}. The analogue of the Ewens Sampling Formula in this case is derived in \citet{GL2005}.

Accuracy of the implementation was checked by ensuring that for smaller sample sizes the equation
\[
p(\bm{n}) = \sum_{s=k-1}^Bp(\bm{n};s)
\]
was approximately satisfied, where $k$ is the number of alleles in $\bm{n}$ and $B$ is a suitable upper bound. The simulation variance, starting with different seeds, was observed to be small for sample sizes such as in example \ref{HammSect2} that follows.

\subsection{Hammer et al. example, continued}\label{HammSect2}

We continue with the example started in Section~\ref{HammSect}. The Y haplotype  data of $n = 1,544$ sequences from \citet{H1998}, had 10 haplotypes and $s = 9$ segregating sites. We continue to use their value of $\theta=2.5$ for illustration.  The 10 haplotype frequencies are 
\begin{center}
\begin{tabular}{cccccccccc}
21& 23& 853& 188& 75 &1 &68& 31 &67& 217
\end{tabular}
\end{center}
in the lineage order shown in Figure 7 of~\citep{H1998}. The average values in the tables below are conditional on the configuration and number of segregating sites, thereby extending the results of Section~\ref{HammSect}.

The Appendix describes the input for the implementation of the method. A command line of
\begin{verbatim}
esf_stl HammerHap.dat 10 9 2.5 1000000 93849 -a
\end{verbatim}
in which the input file \texttt{HammerHap.dat} contains the haplotype frequencies in the order above, produces the  output described below; the  average coalescence times are not shown. 

Two runs with different seeds gave identical output to three significant places, showing some confidence in the output. The probability of obtaining the sample configuration and $s = 9$ segregating sites was $1.4785 \times 10^{-19}$.  As a comparison the probability of the sample configuration, calculated from the Ewens Sampling Formula, was $1.1722 \times 10^{-18}$. The mean TMRCA, conditional on the data,  in coalescent units was $1.15$, which may be compared to the value of $1.21$ obtained in Section~\ref{HammSect}.

\subsubsection{Stationary properties}
Here we record some information about the sample at time $0$.

The conditional expected mutation times in increasing time order are
\begin{center}
\begin{tabular}{cccccccccc}
 0.003& 0.022& 0.039& 0.062& 0.094& 0.142& 0.219& 0.360& 0.675
\end{tabular}
\end{center}
while the conditional expected haplotype  loss times in increasing time order are
\begin{center}
\begin{tabular}{cccccccccc}
 0.003& 0.022& 0.039& 0.062& 0.094& 0.142& 0.219& 0.360& 0.761.
\end{tabular}
\end{center}
 Most of the tree structure has developed by an average time of less than 1.00 coalescent time unit.

The conditional expected haplotype ages, in the order they are listed above, are
\begin{center}
\begin{tabular}{cccccccccc}
 0.051& 0.092& 0.995& 0.406& 0.216& 0.007& 0.201& 0.114& 0.200& 0.446.
\end{tabular}
\end{center}
 These are monotonic in the number of copies of the haplotype in the sample, confirming the intuition that common haplotypes tend to be older.
\subsubsection{Time-varying properties}
At a given time $t$ in the past, the distribution of the configuration, number of mutations, and number of ancestral lineages conditional on the current configuration and number of mutations can be calculated by the importance sampling program. We illustrate this by considering time points $t=0.1,0.5,1.0, 1.5$ and taking averages at those times.

We begin by comparing the conditional distribution of $A_n(t)$ and $S_n(t)$ with the analogous results in Table 1. Additionally, Table 2 shows the average of the number of haplotypes, $K_n(t)$, present at time $t$. 

\begin{table}[h]\label{Table:2}
\begin{center}
\begin{tabular}{|c|c|c|c|}
\hline
$t$ & $K_n(t)$  & $S_n(t)$ & $A_n(t)$ \\
\hline
0 & 10 & 9 &1544 \\
0.1 & 5.09 & 4.09 &19.9\\
0.5 & 1.84 & 0.85 & 3.94\\
1.0 & 0.74 & 0.17 & 1.75\\
1.5 & 0.21 & 0.03 &1.19\\
\hline
\end{tabular}
\end{center}
\caption{Result of 1,000,000 runs for the constant population size coalescent model with a sample of size $n = 1,544$ and $9$ segregating sites, and haplotype frequencies given above. Table shows average values of $K_n(t), A_n(t), S_n(t)$ conditional on the haplotype frequencies and  $S_n = 9$. }
\end{table}

The results in the second and third columns should be compared with those in Table~1; they show qualitatively the same results.  Table 3 shows the relative errors for estimates of $S_n(t)$ and $A_n(t)$. For example, letting superscripts 1 and 2 denote estimates from Tables~1 and~2, the relative error for $S_n(t)$ is $|S_n^1(t)-S_n^2(t)|/\frac{1}{2}(S_n^1(t)+S_n^2(t))$.

\begin{table}[h]
\begin{center}
\begin{tabular}{|c|c|c|}
\hline
$t$ & $S_n(t)$  & $A_n(t)$ \\
\hline
0.1 & 0.25 & 0.02\\
0.5 & 0.75 & 0.06\\
1.0 & 0.75 & 0.02\\
1.5 & 1.2 & 0.05\\
\hline
\end{tabular}
\end{center}
\caption{Relative errors for estimates from Tables~1 and~2}
\end{table}
\begin{table}[htb]
\begin{center}
{\small
\begin{tabular}{|c|cccccccccc|}
\hline
time $t$ & \multicolumn{9}{c}{Haplotype frequency}&\\
\hline
0&
21&23&853&188&75&1&68&31&67&217\\
0.1&
0.227&0.250&11.7&2.30&0.862&0.010&0.777&0.340&0.765&2.69
\cr
0.5&
0.033&0.036&2.63&0.372&0.130&0.002&0.117&0.050&0.115&0.441
\cr
1.0&
0.015&0.011&0.837&0.140&0.050&0.001&0.045&0.020&0.044&0.165
\cr
1.5&
0.006&0.003&0.206&0.043&0.016&0.000&0.014&0.006&0.014&0.051\\
\hline
\end{tabular}
}
\end{center}
\caption{Extant haplotype counts for the Hammer data at different times in the past.}
\end{table}

Table 4 gives the number of haplotype counts at different times $t$ in the past. Haplotypes decrease because of coalescence and types are eventually lost when their defining mutation takes place.
In Table 5  we give the distribution of the number of ancestral lines at different times $t$ in the past, conditional on the current haplotype configuration and the number of segregating sites. The most interesting time configurations are when $t \leq 0.1$; afterwards the number of lineages and number of haplotypes decrease rapidly.
%
%
%
%
\begin{table}[htb]
\begin{center}
{\small
\begin{tabular}{|c|ccccccccc|}
\hline
t&12&13&14&15&16&17&18&19&20\cr
0.1&0.001&0.003&0.009&0.023&0.048&0.083&0.121&0.149&0.159\cr
&21&22&23&24&25&26&27&28&\cr
&0.142&0.111&0.074&0.044&0.022&0.010&0.004&0.001&\cr
\hline
0.5&1&2&3&4&5&6&7&8&\cr
&0.011&0.088&0.259&0.337&0.216&0.073&0.014&0.002&\cr
\hline
1.0&1&2&3&4&5&&&&\cr
&0.426&0.414&0.143&0.016&0.001&&&&\cr
\hline
1.5&1&2&3&&&&&&\cr
&0.826&0.163&0.011&&&&&&\cr
\hline
\end{tabular}
}
\end{center}
\caption{Average number of ancestral lineages in the Hammer data at time $t$ in the past.}
\end{table}

\subsection{A 1000 Genomes Y chromosome dataset}\label{ABE138}

A larger Y chromosome data set comes from the 1000 Genomes Project. An analysis of these data is made in \citet{P2016}, where a phylogeny is constructed. The paper concludes that the data show evidence of expansion. 
As an example of our approach, we consider a subset of this  data set consisting of sequences in the A, B and E haplotype groups. 
These are the three oldest groups in the phylogeny, and are composed of 334 sequences.
There may be explanations other than expansion for the data configuration, such as a non-random choice of individuals, however we will assume a random sample for this example analysis.
For illustration we focus on  the TBL1Y gene, composed of some 180,000bp, and containing 278 biallelic SNPs.  The haplotype configuration, with $\alpha_j$ equal to the number of alleles of multiplicity $j$, is given in Table~\ref{table:st1}.
\begin{table}
\begin{center}
{\small
\begin{tabular}{|ccccccccccc|}
\hline
$\alpha_1$&$\alpha_2$&$\alpha_3$&$\alpha_4$&$\alpha_5$&$\alpha_6$&$\alpha_7$&$\alpha_{14}$&$\alpha_{32}$&$\alpha_{50}$&$\alpha_{61}$
\cr
 107
 &12
 &6
 &1
 &1
 &2
 &1
 &1
 &1
 &1
 &1\\
\hline
\end{tabular}
}
\end{center}
\caption{Allele multiplicities observed in the TBL1Y dataset}\label{table:st1}
\end{table}

Watterson's estimate of $\theta$ \citep{W1975} based on the number of segregating sites $s$ in a constant-sized population is
\[
\widehat{\theta}_W = \frac{s}{\sum_{j=1}^{n-1}1/j} = 44
\]
The maximum likelihood estimate $\widehat{\theta}_E$ based on the Ewens' sampling formula uses $k=\sum_{j=1}^n \alpha_j$, a sufficient statistic for $\theta$. $\widehat{\theta}_E$ satisfies
\[
k=1 + \sum_{j=1}^{n-1}\frac{\widehat{\theta}_E}{\widehat{\theta}_E+j}.
\]
In the TBL1Y dataset $\widehat{\theta}_E=82$. The large number of singletons $\alpha_1$ in the data suggests exponential growth in the population. Growth produces a \emph{star shaped} coalescent tree, which leads to a greater number of singleton sequences.
 The mean number of singletons in the constant size population setting is
\[
\bbbe \left[\alpha_1\right] = \frac{n\theta}{n+\theta-1}.
\]
If $\theta=82$, then $\bbbe \left[\alpha_1\right] = 66$, which is much less than the observed $\alpha_1=107$. Tajima's $D$ \citep{T1983} is given by
\[
D = \frac{\pi - \widehat{\theta}_W}
{\sqrt{\widehat{\text{var}}(\pi - \widehat{\theta}_W)}}
\]
where $\pi$ is the average number of pairwise differences, an unbiased estimate of $\theta$.
 This may be used to test for population growth or other departures from the coalescent model with no growth. Large negative  values of $D$ indicate population growth. In our data $\pi=6.49$ and $D=-2.6$, consistent with expansion. In the Appendix we describe another statistic for testing the no-growth model based on the frequency spectrum, and particularly on the number of singletons,  when $\theta$ is large. \citet{A2009} develops neutrality tests based on the frequency spectrum which generalize tests based on the number of segregating sites and Tajima's $D$.

In this large data set it is difficult to obtain a very precise estimate of $\theta$ and growth rate $\beta$ because there is a large amount of variation in the importance sampling scheme due to the number and length of the sequences. This is not so much an importance sampling issue, but due to the size and structure of the data. Random subsets of the data could be chosen, but the main feature of the data is the large number of singletons. Ancestral inference of the number of ancestral lines at $t$ back would also be difficult to interpret for subsamples. We try a large value $\theta=100$ with different growth rates $\beta$. Growth decreases the variation in the sample, but increases the proportion of singletons, because the coalescent lengths are shortened and the tree is star shaped. Increasing $\theta$ with growth keeps the variation as well as increasing the number of singletons.
The likelihood of the allele configuration and number of segregating sites was calculated for $\theta=100$ with several values of the growth rates $\beta$. Two different replicates each with 10 million runs gave the results in Table~\ref{Table5}.

\begin{table}[h]
\begin{center}
{\small
\begin{tabular}{|cccc|}
\hline
$\beta$&Replicate 1&Replicate 2&Average\cr
\hline
0&		2.0934e-61&   3.4592e-62&1.2197e-61\cr
0.5	&	1.2203e-60&	1.0782e-61&6.6405e-61\cr
1.0	&	4.1350e-60&	1.1497e-60&2.6424e-60\cr
1.5	&	9.2120e-61&	2.8962e-61&6.0541e-61\cr
2.0	&	1.5297e-61&   5.9120e-61&3.7209e-61\cr
2.5	&	4.9886e-62&	3.7215e-62&4.3551e-62\\
\hline
\end{tabular}
}
\end{center}
\caption{Likelihoods with expansion in the TBL1Y data.}\label{Table5}
\end{table}

A plausible maximum likelihood estimate when $\theta=100$ is $\widehat{\beta }=1.0$. If $\theta$ and $\beta$ are increased together, it is possible that the likelihood estimates fall on a ridge. In the first replicate with these values of $\theta$ and $\beta$ the TMRCA was $1.461$ and the average of ages within halpotype groups are shown in Table~\ref{Table6}.
\begin{table}[h]
\begin{center}
{\small
\begin{tabular}{|cccccc|}
\hline
& &\multicolumn{2}{c}{Haplotype groups} &&\\
\hline
$\alpha_1$ & $\alpha_2$ & $\alpha_3$ & $\alpha_4$ & $\alpha_5$ & $\alpha_6$\\
0.0123 & 0.0086 & 0.0222 & 0.0249 & 0.0421& 0.0390 \\
\hline
$\alpha_7$&$\alpha_{14}$&$\alpha_{32}$&$\alpha_{50}$&$\alpha_{61}$ &\\
0.0290&0.1338& 0.1331& 0.1259& 0.1304 &\\
\hline
\end{tabular}
}
\end{center}
\caption{Allele age within groups.}\label{Table6}
\end{table}

With the large value of $\theta=100$, times where mutation creates an allele are close to the leaves of the tree.

The scale of modern molecular datasets points out the difficulty of \emph{exact} inference techniques, and highlights the need for alternative approaches. Among these are the Approximate Bayesian Computation (ABC) and other advanced simulation methods, to which Paul Joyce made several contributions. Two examples of Paul's research are in estimating evolutionary rates of trait evolution by ABC \citep{SHJA2012}, and a perfect simulation method for simulation from a non-neutral high dimensional distribution of allele frequencies \citep{JGB2012}.
 
\section{Acknowledgements}

We thank Chris Tyler-Smith,  Yali Xue and two referees for their helpful comments about the 1000 Genomes Y chromosome data. 

Code for the rejection and importance sampling methods may be obtained from the authors. 

\section*{Appendix}

\subsection*{Importance sampling code}

The input for the importance sampling method illustrated in Section~\ref{HammSect2} is:
\begin{verbatim}
esf_stl configfile k [#alleles] m [#mutations] theta replicates seed

	Options
	-g beta [exponential growth]
	-a [age information]
	-t time [Configuration at time]
	
	Bob Griffiths	4 May 2017, Version 1.7
\end{verbatim}

\subsection*{Poisson approximation for large $\theta$}

Motivated by the discussion in Section~\ref{ABE138} we discuss the behaviour of the Ewens sampling Formula for large values of $\theta$ and $n$.
The Ewens Sampling Formula gives the distribution of the number of haplotypes and their frequencies in a sample taken from a constant-size population. Writing $\alpha_j$  for the number of haplotypes with frequency $j$,  the distribution is 
\begin{equation}\label{ESF1}
p_E(\alpha_1,\alpha_2,\ldots,\alpha_n) = \frac{n!}{\theta_{(n)}} \, \prod_{j=1}^n \left(\frac{\theta}{j}\right)^{\alpha_j} \, \frac{1}{\alpha_j!},
\end{equation}
where $\theta_{(n)} := \theta(\theta+1) \cdots (\theta+n-1)$, and $\alpha_1 + 2 \alpha_2 + \cdots + n \alpha_n = n$.
The formula (\ref{ESF1}) shows that, were it not for the condition that $\alpha_1 + 2 \alpha_2 + \cdots + n \alpha_n = n$, the $\alpha_j$ would be independent Poisson random variables with mean $\theta/j$. Indeed, for fixed $\theta$ it is known that for any $b = o(n)$ as $n \to \infty$, the total variation distance between the distribution of $(\alpha_1,\ldots,\alpha_b)$ and that of $(Z_1,\ldots,Z_b)$, for independent Poisson random variables with $\bbbe \left[Z_j\right] = \theta/j$, is $O(b/n)$ as $n \to \infty$. See~\citet{ABT}, Theorem 5.2.

Here we consider the case in which $\theta \to \infty$ with $n$, and we  show that for fixed $b$, $(\alpha_1,\ldots,\alpha_b)$  has asymptotically  the distribution of $(Z_1^\theta,\ldots,Z_b^\theta)$, where
$$
\bbbe Z_j^\theta = \frac{\theta}{j}\, \left( \frac{n}{n+\theta}\right)^j, \quad j=1,\ldots,b.
$$

To see this, consider the joint falling factorial moments of $\alpha_1,\ldots,\alpha_b$, given by~\citet{W1974} as
\begin{eqnarray*}
\bbbe \prod_{j = 1} ^ b (\alpha_j)_{[r_j]} & = & \mathbbm{1}(m \leq n)\, \frac{n!\, \Gamma(\theta+n-m)}{(n-m)!\, \Gamma(\theta+n)} \, \prod_{j=1}^b\left( \frac{\theta}{j}\right)^{r_j},\\
& = & \mathbbm{1}(m \leq n)\, \frac{n(n-1)\cdots(n-m)}{(\theta+n-1) \cdots(\theta + n - m)} \, \prod_{j=1}^b\left( \frac{\theta}{j}\right)^{r_j}\\
& \sim & \left(\frac{\theta}{\theta+n}\right)^m\, \prod_{j=1}^b\left( \frac{\theta}{j}\right)^{r_j}\\
& = & \prod_{j=1}^b\left( \frac{\theta}{j}\left( \frac{n}{\theta+n}\right)\right)^{r_j}\
\end{eqnarray*}
where $m := r_1 + 2 r_2 + \cdots + b r_b$. 
The term on the right gives the falling factorial moments of $(Z_1^\theta,\ldots,Z_b^\theta)$, and the result follows from the method of moments.
 
 
In practice, different limit laws are obtained depending on the way $\theta$ varies with $n$. For example, if $\theta \sim \eta n$, then $\alpha_1$ has approximately a Poisson distribution with mean $\theta/(1+\eta)$. For the data in Section~\ref{ABE138}, with $\theta = 82$, the number of singleton haplotypes has mean $82 \cdot (334/416) \approx 65.84$. Since the probability of observing 107 or more singletons is then $\approx 1.92 \times 10^{-6}$, we conclude that the constant-size model does not provide an adequate fit.
In a similar spirit, $\alpha_1 + \alpha_2$ has approximately a Poisson distribution with mean 92.27. We observed $\alpha_1 + \alpha_2 = 119,$ the probability of a larger value being $\approx 0.0043$; once more, this suggests the constant-size model is not a good fit. 

\section{References}

\end{document}